\documentclass{article}
\usepackage{amsfonts}
\usepackage{amsmath}
\usepackage{amssymb}
\usepackage{graphicx}
\usepackage{rotating}
\usepackage{lscape}
\usepackage{longtable}
\usepackage{epsfig}
\usepackage[all]{xy}
\newcounter{excnt}
\newtheorem{thm}{Theorem}
\newtheorem{df}[thm]{Definition}
\newtheorem{rem}[thm]{Remark}

     {\noindent {\bf Proof:}}{\nopagebreak\begin{flushright}{\bf q.e.d.}\end{flushright}}
\newenvironment{example}%
     {\stepcounter{excnt}
      \noindent {\bf Example \arabic{excnt}} \phantom{i}}{\vspace{0.4cm}}

\title{A Short Note on Hauser's Kangaroo Phenomena and 
       Weak Maximal Contact in Higher Dimensions}
\author{Anne Fr\"uhbis--Kr\"uger}

\begin{document}

\maketitle

\begin{abstract}
Currently there are several approaches to resolution of singularities in 
positive characteristic all of which have hit some obstruction.
One natural idea is to try to construct new meaningful examples at this
point to gain a wider range of experience. To produce such examples
we mimic the characteristic zero approach and focus on cases where it fails.
In particular, this short note deals with an example-driven study of
failure of maximal contact and the search for an appropriate replacement.
\end{abstract}

\section{Introduction}

Hypersurfaces of maximal contact are one of the key concepts in Hironaka's
inductive proof of desingularization in characteristic zero, but unfortunately
they need not even exist locally in positive characteristic as e.g. 
Narasimhan's example \cite{Nar} shows. 
In \cite{EH} and \cite{Ha2} 
Hauser replaces hypersurfaces of maximal contact by the characteristic-free
notion of hypersurfaces of weak maximal contact, i.e. hypersurfaces which 
maximize the order of the subsequent coefficient ideal, but which do not 
necessarily contain the equiconstant points after all sequences of blowing 
ups in permissible centers. In the corresponding approach (\cite{Ha},
\cite{HW}) 
to resolution of surface singularities in positive characteristic, 
this modification of the concept of maximal contact turns out to be sufficient
to enter into an approach in the flavour of Hironaka's original induction on 
the dimension of the ambient space. To obtain desingularisation 
of surfaces along those lines, this is, of course, not the only change to the 
characteristic zero arguments; important further modifications to certain 
components of the desingularisation invariant are required. Considering higher
dimensions, however, the first step toward a construction of a 
desingularization similar to the characteristic zero approach or even toward 
new meaningful examples illustrating the obstructions against it again needs 
to be a reconsideration of the right generalization of maximal contact.\\
\\
For readers convenience, we briefly recall some key concepts in section 2. 
Here one focus will be on the question of recognition of a potential kangaroo.
In section 3, we start by considering an example where the original 
definition of weak maximal contact does not suffice for the description of 
a kangaroo phenomenon and then suggest a slightly modified version which 
is suitable for any dimension and not just surfaces. Using this new notion 
of a flag of weak maximal contact, section 4 is then devoted to examples of 
the different roles which the hypersurfaces originating from the flag 
can play in the course of a sequence of permissible blowing ups.\\
\\
The author would like to thank Herwig Hauser, Vincent Cossart, 
Dominique Wagner and Santiago Encinas for fruitful discussions which
originally started her interest in desingularisation in positive 
characteristic. The author is also indebted to the referee of the article
whose comments were a great help in revising a previous version 
of this article. All examples appearing in this article originated from
structured experiments using the computer algebra system {\sc Singular}. 

\section{Basic facts and definitions}

A section of just a few pages is obviously not sufficient to even
give a brief overview of the tools and general philosophy of algorithmic 
desingularization, let alone all the delicacies of the case of positive 
characteristic. On the other hand, more than just 5 pages would be by far
too long compared to the following two sections. Hence we do not attempt 
this here, but only very briefly sketch the idea of the characteristic 
zero resolution process to give a context, subsequently recalling the 
notions of hypersurfaces of weak maximal contact and of kangaroo points 
in positive characteristic. For additional background information on 
the characteristic zero case, we would like to point to more thorough 
discussions in section 4.2 of \cite{FK1} 
from the practical point of view and in \cite{EH} embedded in a detailed
treatement of the resolution process. For a detailed introduction to 
characteristic $p$ phenomena and kangaroo points see \cite{Ha}.

\subsection{The philosophy of the characteristic zero approach}

In Hironaka's original work \cite{Hir} and in all algorithmic approaches 
based on it, e.g. \cite{BM},\cite{BEV},\cite{EH}, the general approach
is that of a finite sequence of blow-ups at appropriate non-singular
centers. The very heart of these proofs is the choice of center which is
controlled by a tuple of invariants assigned to each point; it is of a 
structure similar\footnote{In the case of Bierstone and Milman, the very
first entry is a finer invariant, the Hilbert-Samuel function.} 
to the following one
               $$(ord,n;ord,n;\dots)$$
with lexicographic comparison, the upcoming center being the set of maximal
value of the invariant. Here $ord$ stands for an order of an appropriate 
(auxiliary) ideal (see below), $n$ for a counting of certain exceptional 
divisors. At each ';' a new auxiliary ideal of smaller ambient dimension, 
a coefficient ideal, is created by means of a hypersurface of maximal 
contact.\\

To fix notation, let $W$ be a smooth equidimensional
scheme over an algebraically closed field $K$ of characteristic zero and
$X \subset W$ a subscheme thereof. We now immediately focus on one 
affine chart $U$ with coordinate ring $R$ and denote the maximal ideal 
at $x \in U$ by ${\mathfrak m}_x$.\\
The order of the ideal $I_X = \langle g_1,\dots,g_r \rangle \subset R$ at a
point $x \in U$ is defined as
           $$ord_x(I):={\rm max}\{m \in {\mathbb N} | \;\;
                   I \subset {\mathfrak m}_x^m \}.$$
In characteristic zero, the order of the non-monomial part of an ideal
can never increase under blow-ups which makes it a good ingredient
for the controlling invariant of the resolution process whose decrease
marks the improvement of the singularities.\\

For the descent in ambient dimension, hypersurfaces of maximal contact
are required; these locally contain all points of maximal order,
satisfy certain normal crossing conditions and continue to contain all
points at which the maximal order did not yet drop after 
any permissible sequence of blow-ups. In characteristic zero, they
always exist locally and can be computed in a rather straight-forward
way. The construction of the coefficient ideal for $I$ at $X$ w.r.t. 
a hypersurface of maximal contact $Z=V(z)$ is then performed in the
following way:
     $$Coeff_Z(I)=\sum_{k=0}^{ord_x(I)-1} I_k^{\frac{k!}{k-i}}$$
where $I_k$ is the ideal generated by all polynomials which appear as
coefficients of $z^k$ in some element of $I$. Given this notion of 
coefficient ideal, it is possible to rephrase the condition on a 
hypersurface of maximal contact from 'containing all points of
maximal order' to 'maximizing the order of the non-monomial part of
the arising coefficient ideal under all choices of hypersurfaces'.

\subsection{Weak maximal contact and kangaroos}

In positive characteristic, there are well known examples of failure of 
maximal contact in the sense that eventually the equiconstant points
will leave the strict transform of any chosen smooth hypersurface
(see \cite{Nar}). Using the characteristic free formulation of the
first condition for maximal contact, i.e. that it should maximize the order
of the non-monomial part of the subsequent coefficient ideal, and dropping 
the condition that this should hold after any permissible sequence of 
blow-ups, we obtain Hauser's definition of weak maximal contact. In this way,
Hauser and Wagner \cite{HW} then allow passage to a new hypersurface of
weak maximal contact, if the previously chosen one happens to fail to 
have the maximizing property at some moment in the resolution process.\\

Additionally there are examples (see \cite{Moh}) in which the order of
the non-monomial part of the first coefficient ideal can increase under 
a sequence of blow-ups in positive characteristic. In \cite{Ha}
Hauser shows that these two phenomena are closely related in the sense
that both arise in the same rather rare settings and gives an
explicit criterion for the possibility of such a phenomenon, 
which he calls a kangaroo point focusing on the point where this occurs.
In this article, we often choose to refer to this as a kangaroo phenomenon,
emphasizing the fact that not the point itself is in the center of interest,
but the deviation from the characteristic zero case. Using the same notation
for $W$, $X$ etc. as in the previous section, we now recall Hauser's 
definition: 

\begin{df}[\cite{Ha}]
Let $\pi: W' \longrightarrow W$ be a blow-up at a permissible center $Z$,
and $x \in Z$ a point of maximal order $c$ for $I_X$. Denoting the weak
transform of $X$ under $\pi$ by $X'$, let $x'\in X' \cap \pi^{-1}(x)$ be a 
point at which $ord_{x'}(I_{X'})=c$. Then $x'$ is called a kangaroo point, 
if the order of the non-monomial part of the coefficient ideal of 
$I_X$ at $x$ w.r.t. a hypersurface of weak maximal contact is less than 
the order of the non-monomial part of the coefficient ideal of $I_{X'}$ 
w.r.t. a (possibly newly chosen) hypersurface of weak maximal contact.
\end{df}

\begin{df}
Generalizing Hauser's notion of a kangaroo point, we shall call a blowing
up, at which such an increase in order occurs for one of the coefficient
ideals at some level in the descent of ambient dimension, a kangaroo
phenomenon.
\end{df}

\begin{rem}[\cite{Ha}] 
A kangaroo point can only occur, if the following conditions are satisfied:
\begin{itemize}
\item[(a)] the order $c$ of the ideal $I_{X}$ at $x$ does not exceed 
           the order of $I_{X'}$ at $x'$
           and is divisible by the charateristic of the ground field.
\item[(b)] The order of the non-monomial part of the coefficient 
           ideal is a multiple of $c$.
           \footnote{For kangaroo phenomena, this condition should 
           analogously read
           'one of the coefficient ideals occurring in the descent of 
            ambient dimension'.} 
\item[(c)] The exceptional multiplicities of the coefficient 
           ideal need to satisfy a certain numerical inequality 
           (whose specification would need to much room here).
\end{itemize}
\end{rem}

This remark does not yield a sufficient criterion of detection of kangaroos. 
However, if a kangaroo phenomenon occurs, then its effect is an increase 
of order of the non-monomial part of the coefficient ideal by means of 
leaving at least two exceptional divisors at the same time and a suitable 
change of hypersurface of weak maximal contact (see examples in sections 
3 and 4 for details).\\

Combining the above observations of Hauser with well-known observations
by Hironaka and Giraud, condition (a) can be made a bit more precise. 
To this end, we need to recall another singularity invariant, the ridge 
(french: la fa\^ite). Following the exposition of \cite{Oda}, let us 
consider the tangent cone $C_{X,x}$ of $I_{X}$ at $x$ and the 
largest subgroup scheme $A_{X,x}$ of the tangent space $T_{W,x}$
satisfying the conditions that it is homogeneous and leaves the 
tangent cone stable w.r.t. the translation action. $A_{X,x}$ is called the 
ridge of the tangent cone of $I_X$ at $x$.\\

It is a well-known, important fact that the ridge can be generated by 
additive polynomials, i.e. by polynomials of the form
$$ \sum_{i=1}^n a_ix_i^{p^e}$$
where $p$ is the characteristic of the underlying field. In characteristic 
zero the ridge is always generated by polynomials of degree one; in 
positive characteristic the occurrence of a ridge not generated by polynomials
of degree one marks a point for which the reasoning of characteristic zero 
might break down. Following the exposition of \cite{BHM} 
the ridge can also be phrased as the smallest set of additive polynomials
$\{p_1,\dots,p_r\}$ generating the smallest algebra $k[p_1,\dots,p_r]$
such that 
$$I_X = (I_X \cap k[p_1,\dots,p_r])k[\underline{x}].$$\\

Combining this with Hauser's condition (a), we obtain a refined version for
hypersurfaces, which, of course, still requires 
$ord_{x'} (I_{X'}) = c = ord_x (I_X)$ 
and, additionally, that the ridge must at least have one generator in higher 
degree, i.e. in some degree $p^e$. This sharpens the condition of divisibility
of the order by a p-th power to the fact that some variable actually 
only occurs as p-th powers in the tangent cone and is implicitly already 
present in \cite{Ha}. 
According to Hauser's condition (b), the degree of the non-monomial part of
the first coefficient ideal is required to be a multiple of the degree $c$.
In contrast to condition (a), this can not be made more precise by simply
adding the condition that the ridge of the non-monomial part of this 
coefficient ideal is not generated in degree 1, because higher order generators
of the coefficient ideal might introduce lower degree polynomials into the ridge
which allow dropping of certain contributions arising from the lowest order generators
of the ideal. To illustrate the role of the ridge, we give 3 examples:\\

\begin{example}
Over a field $K$ of characteristic 3, consider an affine chart 
$U={\mathbb A}_K^4$ (with variables named $x,y,z,w$) which already 
results from a sequence of 2 blow-ups and contains exceptional divisors 
$E_1=V(w)$ and $E_2=V(z)$, born from the first and second blow-up respectively.
(These two blow-ups are indeed necessary for the possibility of an 
occurrence of a kangaroo point after the subsequent blowing up, according 
to Hauser's technical condition (c) which was not formulated explicitly 
in the previously stated remark.) \\ 
Locally at the coordinate origin of this chart, consider the three 
subvarieties of ${\mathbb A}_K^4$ defined by the following ideals:
\begin{itemize}
\item $I_{X_1}=\langle x^3+z^{14}w^{10}(z^6-w^6) \rangle$ \\
      This is the strict transform\footnote{Actually this is the weak 
      transform of $I_{X_1}$ which in the principal ideal case happens
      to coincide with the strict transform.} 
      of $\langle x^3+z^{13}-zw^{18} \rangle$ under the two blow-ups. 
      The ridge of $I_{X_1}$ can obviously be described by $\{x^3\}$, the 
      non-monomial part of its first coefficient ideal is
      $$\langle z^{12}+z^6w^6+w^{12} \rangle,$$
      with ridge $\{z^3,w^3 \}$. \\
      After blowing up again at the
      origin, we obtain (in the $E_3=V(w)$-chart) the strict transform 
      $$I_{X_1'}= \langle x^3+z^{14}w^{27}(z^6-1) \rangle$$
      which after a coordinate change $z_{new}=z-1$ and a passage to
      a new hypersurface of weak maximal contact 
      $V(x+z_{new}^2w^9)=V(x_{new})$ reads as
      $I_{transf.}= 
              \langle x_{new}^3+z_{new}^6w^{27}(-z_{new}+h.o.t.) \rangle$.
      Since $z_{new}$ does not correspond to an exceptional divisor, this has 
      a non-monomial part of the first coefficient ideal of the form
      $$\langle z_{new}^{14}+h.o.t.\rangle.$$
      This ideal is of order 14 as compared to the corresponding order 12 
      before the last blowing up which clearly indicates the occurrence 
      of a kangaroo point.
\item $I_{X_2}=\langle x^2y+z^{14}w^{10}(z^6-w^6) \rangle$ \\
      This is the strict transform 
      of $\langle x^2y+z^{13}-zw^{18} \rangle$ under the two 
      blow-ups.\footnote{Here we are actually already deviating a bit from 
      Hauser's original definition, because we consider an initial part
      involving 2 variables and then descend in ambient dimension in one
      step of 2 to $V(x,y)$ seen as a hypersurface in $V(x)$ which is
      in turn a hypersurface in ${\mathbb A}^4$. This is possible by
      collecting all coefficients of monomials of the form $x^ay^b$ with
      $a+b=k$ into the ideal $I_k$; for more details on this see e.g.
      \cite{FK2}, where this has been used in a very explicit way.} 
      The ridge of $I_{X_2}$ is obviously $\{x,y \}$, the 
      non-monomial part of its coefficient ideal w.r.t. the descent
      in ambient dimension to $V(x,y)$ is
      $$\langle z^{12}+z^6w^6+w^{12} \rangle$$
      as before with ridge $\{ z^3,w^3 \}$. \\
      After blowing up again at the origin, we obtain 
      (in the $E_3=V(w)$-chart) the strict transform 
      $$I_{X_2'}= \langle x^2y+z^{14}w^{27}(z^6-1) \rangle$$
      for which even a coordinate change $z_{new}=z-1$ cannot lead to a 
      kangaroo point, because no suitable passage to new hypersurfaces 
      of weak maximal contact killing the term 
      $z_{new}^6w^{27}$ is available. 
      This could already be expected at the beginning due to the fact that the
      ridge of $I_{X_2}$ is generated in degree 1.
\end{itemize}

The third example is of a different flavor and only serves to illustrate,
how higher order generators of the ideal might influence the ridge
in a way which is not desirable for the consideration of coefficient ideals:
\begin{itemize}
\item $I_{X_3}=\langle x^3+z^{14}w^{10}(z^6-w^6), 
                       z^{30}w^{17}(y^{19}+y^5z^7w^3)\rangle$\\
      This is the weak transform of $\langle x^3+z^{13}-zw^{18}, 
      y^5z^{18}+y^{19}w \rangle$ under the two blow-ups.
      The ridge of $I_{X_3}$ is obviously $\{ x^3,y,z,w \}$, whereas
      only the hypersurface $V(x)$ can be chosen as hypersurface of weak
      maximal contact. The non-monomial part of the first coefficient ideal is
      \begin{center} $\langle z^{12}+z^6w^6+w^{12}, 
        (z^6-w^6)(z^{16}w^{7}(y^{19}+y^5z^7w^3)),$\\
        $ z^{32}w^{14}(y^{38}-y^{24}z^7w^3+y^{10}z^{14}w^6) \rangle.$ 
      \end{center}
      The ridge can be computed to be $\{ y,z,w \}$, e.g. by the
      algorithm of \cite{BHM}. \\
      After blowing up again at the
      origin, we obtain (in the $E_3=V(w)$-chart) the weak transform 
      $$I_{X_3'}= \langle x^3+z^{14}w^{27}(z^6-1),\dots \rangle$$
      which after a coordinate change $z_{new}=z-1$ and a passage to
      a new hypersurface of weak maximal contact 
      $V(x+z_{new}^2w^9)=V(x_{new})$ has the same first generator of order $14$ 
      as in example 1, the second generator does not have effect on the 
      order of the non-monomial part of the first coefficient ideal as
      can be checked by explicit computation. Comparing this to the
      first example, we see that the higher order generator, which 
      does not actually influence the order of the non-monomial part of
      the coefficient ideal, masked the situation in the computation of
      the ridge.
\end{itemize}  
From these three examples, we see the usefulness of the ridge for 
anticipating kangaroo points in the case of hypersurfaces, whereas 
in the case of ideals this may be hidden by contributions of higher order 
generators. However, if we only consider the ridge of the ideal which
is generated precisely by the lowest-order generators of the original ideal
(instead of the ridge of the whole ideal), then there is hope to use this 
new ridge for ideals and maybe even to slightly sharpen item (b) 
in Hauser's condition for kangaroo points.  
\end{example}

\begin{rem}
These considerations already suggest a strategy for finding interesting
examples by constructing hypersurfaces for which the ridge is not
generated in degree 1 and, additionally, at least once during the iterated 
descents in ambient dimension the ridge of the ideal generated by the 
lowest order generators (denoted from now on as n-ridge for short) 
of the non-monomial part of the respective coefficient ideal is also not 
generated in degree one. In the experiments, which lead to the examples 
of the subsequent sections, an additional heuristic in the choice of 
hypersurfaces of weak maximal contact was used: 
When given the choice between different hypersurfaces, more precisely 
between linearly independent initial parts of possible hypersurfaces, 
we try to minimize the degree of the generator of the ridge/n-ridge 
corresponding to the chosen hypersurface. The reasoning behind this 
heuristic is to force the unpleasant, but interesting behaviour into the 
lowest possible ambient dimension and hence keep a clearer view of the 
occurring phenomena. 
\end{rem}

\begin{rem}
Similar examples to those of the subsequent sections can easily be constructed
in any positive characteristic. For section 3 this is straight forward,
for section 4 it is best achieved by starting in the middle, i.e. precisely
where the first kangaroo has just occurred and construct from there by blowing
down and blowing up.
\end{rem}

\section{In higher dimension not all 
                    hypersurfaces of weak maximal contact are suitable}

The following example shows that the property of maximizing the order of the 
non-monomial part of the upcoming coefficient ideal is not sufficient to 
properly cover all kangaroo phenomena in higher dimensions. It is stated 
in characteristic 2 to allow considerations in rather low degrees, but 
similar examples can be constructed for any positive characteristic.\\

\begin{example}
We consider a sequence of three blow ups of the hypersurface
$V(x^2+w^3+y^{25}+yz^{16}) \subset {\mathbb A}_{K}^4$, $char(K)=2$, 
$K=\overline{K}$.
At each step the respective maximal orders, chosen 
hypersurfaces of weak maximal contact and coefficient ideals are specified. 
In the presence of exceptional divisors, we make use of Bodnar's trick 
\cite{Bod}, which allows skipping the intersection with exceptional 
divisors in intermediate levels of the descent in ambient dimension, if 
we have normal crossing between the upcoming hypersurface of weak maximal 
contact and the exceptional divisors. \\[0.5cm]
To keep the whole rather lengthy sequence of blowing ups more readable, we only
give rather scarce comments. A more commented version of a single blowing up
step was already stated at the end of the previous section.\\
\\
\noindent
{\bf original hypersurface}: \\
     $I=\langle f \rangle = \langle x^2+w^3+y^{25}+yz^{16} \rangle$
\begin{itemize}
\item in ambient space ${\mathbb A}_K^4$\\
      $I=\langle x^2+w^3+y^{25}+yz^{16} \rangle$\\
      The maximal order 2 is attained at $V(x,y,z,w)$.\\
      The ridge of this ideal corresponds to $\{x^2\}$.\\
      As hypersurface of weak maximal contact we may use $H_1=V(x)\subset{\mathbb A}_K^4$.
\item in ambient space $H_1$\\
      $I_{H_1}=\langle w^3+y^{25}+yz^{16} \rangle$.\\
      The maximal order of 3 is then again attained at the origin of $H_1$.\\
      The n-ridge (in the short-hand notation introduced in section 2) is $\{w\}$\\  
      As hypersurface of weak maximal contact we now use \\
      $H_2=V(x,w)\subset H_1 \subset {\mathbb A}_K^4$.
\item in ambient space $H_2$\\
      $I_{H_2}=\langle y^{50}+y^2z^{32} \rangle = \langle (y^{25}+yz^{16})^2 \rangle$.\\
      The maximal order of 34 is again attained at the origin of $H_2$
      and the n-ridge is $\{y^2,z^{32}\}$. 
\item The only possible choice of center is $V(x,y,z,w)$.
\end{itemize}
      As a sideremark to the coefficient ideal in ambient space $H_2$: Here it becomes evident 
      that there are 2 mechanisms which can cause the n-ridge to have generators in higher 
      degree: on one hand, it may be an honest generator in higher degree, on the 
      other hand, it might have arisen from taking powers of contributing ideals $I_k$
      when forming the coefficient ideal (see section 2). However, taking powers can not 
      accidentally cause the degree of a generator of the ridge to drop. Hence the degree of the 
      generators of the ridge can still be used as a rather weak indicator for the 
      possibility of new phenomena in characteristic $p$. Moreover, a higher degree  
      generator of the n-ridge arising from mechanism 2 is only likely to occur, if the 
      contributing ideals $I_k$ are principal, because otherwise mixed products
      of generators would exist in the set of generators of the power of $I_k$. \\ 

\noindent
{\bf after first blowing up, chart $E_1=V(y)$}:\\ 
     $I_{strict}= \langle x^2+yw^3+y^{23}+y^{15}z^{16} \rangle$
\begin{itemize}
\item in ambient space ${\mathbb A}_K^4$\\
      $I_{strict}= \langle x^2+y(w^3+y^{22}+y^{14}z^{16}) \rangle$\\
      The maximal order is again $2$, attained at the origin and the ridge is 
      again $\{x^2\}$. We can keep the strict transform of $H_1$ as our 
      hypersurface of weak maximal contact. (As $\{E_1,{H_1}_{strict}\}$ has
      normal crossings, we may use Bodn\`ar's trick \cite{Bod} and hand the 
      exceptional divisor down to the lower dimension instead of intersecting
      with it at this point.)
\item in ambient space ${H_1}_{strict}$\\
      The non-monomial part of the coefficient ideal\footnote{As taking the coefficient
      ideal and subsequently calculating the controlled transform under the blowing up 
      on one hand and calculating the weak transform of the ideal followed by computing
      the new coefficient ideal on the other hand are known (e.g. \cite{EH}) to lead 
      to the same ideal, we won't go into details on this point.} is
      $\langle w^3+y^{22}+y^{14}z^{16} \rangle$.\\
      The maximal order of $3$ is again attained at the origin and the n-ridge is
      $\{w\}$ as before. We may also use the strict transform of $H_2$ again for the
      descent in ambient dimension. (Here we have normal crossing of $\{E_1,{H_1}_{strict},
      {H_2}_{strict}\}$ and can again use Bodn\`ar's trick.)
\item in ambient space ${H_2}_{strict}$\\
      non-monomial part of coefficient ideal: 
               $\langle y^{16}+z^{32} \rangle = \langle (y^8+z^{16})^2 \rangle$\\
      maximal order $16$ attained at the origin\\
      n-ridge: $\{y^{16}\}$
\item It is easy to check that here again the choice of center has to be the origin.
\end{itemize}

\noindent
{\bf after second blowing up, chart $E_2=V(z)$}:\\
      $I_{strict}= \langle x^2+yz^2w^3+y^{23}z^{21}+y^{15}z^{29} \rangle$
\begin{itemize}
\item in ambient space ${\mathbb A}_K^4$\\
      $I_{strict}= \langle x^2+yz^2(w^3+y^{22}z^{19}+y^{14}z^{27}) \rangle$\\
      maximal order: $2$ at $V(x,zw,yz)$\\
      ridge: $\{x^2\}$\\
      hypersurface of weak maximal contact: strict transform of $H_1$\\
      ($\{{E_1}_{strict},E_2,{H_1}_{strict}\}$ n.cr.)
\item in ambient space ${H_1}_{strict}$\\
      non-monomial part of coefficient ideal: 
               $\langle w^3+y^{14}z^{19}(y^8+z^{16}) \rangle$\\
      maximal order: $3$ at $V(w,yz)$
      n-ridge: $\{w^3\}$
      hypersurface of weak maximal contact: strict transform of $H_2$\\
      ($\{{E_1}_{strict},E_2,{H_2}_{strict}\}$ n.cr.)
\item in ambient space ${H_2}_{strict}$\\
      non-monomial part of coefficient ideal:
               $\langle y^{16}+z^{16} \rangle = \langle (y^8+z^8)^2 \rangle$\\
      maximal order: $16$ at $V(y+z)$\\
      n-ridge: $\{y^{16}+z^{16}\}$
\item center needs to be $V(x,y,z,w)$ as the locus of maximal order after the second
      descent in ambient dimension is not normal crossing with the exceptional divisors
\end{itemize}

\noindent
{\bf after third blowing up, chart $E_3=V(z)$}:\\
      $I_{strict}=\langle x^2+yz^4w^3+y^{23}z^{42}+y^{15}z^{42} \rangle$
\begin{itemize}
\item in ambient space ${\mathbb A}_K^4$\\
      $I_{strict}=\langle x^2+yz^4(w^3+y^{22}z^{38}+y^{14}z^{38}) \rangle$\\
      maximal order: $2$ at $V(x,zw)$\\
      ridge: $\{x^2\}$\\
      hypersurface of weak maximal contact: strict transform of $H_1$\\
      ($E_1$ does not meet this chart, $\{{E_2}_{strict},E_3,{H_1}_{strict}\}$ n.cr.)
\item in ambient space ${H_1}_{strict}$\\
      non-monomial part of coefficient ideal:
               $\langle w^3+y^{12}z^{42}(y^8+1) \rangle$\\
      maximal order: $3$ at $V(w,yz(y+1))$
      n-ridge: $\{w\}$\\
      hypersurface of weak maximal contact: strict transform of $H_2$\\
      ($\{{E_1}_{strict},E_2,{H_2}_{strict}\}$ n.cr.)
\item in ambient space ${H_2}_{strict}$\\
      non-monomial part of coefficient ideal:
               $\langle y^{16}+1 \rangle = \langle (y^8+1)^2 \rangle$\\
      maximal order: $16$ at $V(y+1)$ 
\end{itemize}
Changing the hypersurface for the first descent in ambient dimension from ${H_1}_{strict}$
to $V(x+(y+1)^4z^{21})$, however, we may increase the order of the coefficient ideal in 
ambient dimension $2$. For simplicity of notation, we first make a coordinate change which
translates the point of maximal order to the coordinate origin:
\begin{itemize}
\item in ambient space ${\mathbb A}_K^4$\\
      $\langle x^2+z^4(w^3+ y_{new}w^3 + y_{new}^8z^{38} + y_{new}^9z^{38} 
                  +h.o.t.) \rangle$\\
      maximal order $2$ at $V(x,zw)$ \\
      ridge: $\{x^2\}$\\
      new hypersurface of weak maximal contact: ${H'_1}=V(x+y^4z^{21})$\\
      ($E_1$ does not meet this chart, $\{{E_2}_{strict},E_3,{H'_1}\}$ n.cr.)
\item in ambient space ${H'_1}$\\
      non-monomial part of coefficient ideal:
           $\langle w^3 + y_{new}^9z^{38}+ h.o.t. \rangle$\\
      maximal order $3$ at $V(w,y_{new}z)$  \\
      n-ridge: $V(w)$\\
      hypersurface of weak maximal contact: ${H'_2}=V(w)$\\
      ($\{{E_2}_{strict},E_3,{H'_2}\}$ n.cr.)
\item in ambient space $H'_2$\\
      non-monomial part of the coefficient ideal:
               $\langle y_{new}^{18} + h.o.t. \rangle$\\
      maximal order: $18$ {\bf exceeds previous order $16$}\\
      {\bf kangaroo} phenomenon
\end{itemize}

Here the new phenomenon is that the change of the hypersurface of weak maximal
contact was not forced by the first coefficient ideal, but by one of the later
ones which would not be covered by the standard definition of weak maximal 
contact. 
\end{example}

In the light of the previous example, we suggest a slightly modified version
of weak maximal contact:

\begin{df}
Consider a given point $x$ of a scheme $X$ (possibly in the presence of an 
exceptional divisor $E$) and pass to an affine chart $U$ containing this 
point. We call a flag 
$$ {\mathcal H} = H_1 \supset H_2 \supset \dots \supset H_s $$
admissible at $x$, if the following properties hold:
\begin{itemize}
\item[(a)] $H_1$ is a smooth hypersurface in the ambient space $U$. 
           $H_{i+1}$ is a smooth hypersurface in $H_i$.
\item[(b)] $H_i$ is a hypersurface of weak maximal contact for the coefficient
           ideal obtained by descent of the ambient space through
           $H_1,\dots,H_{i-1}$.
\item[(c)] $x \in H_s$.
\end{itemize}
${\mathcal H}$ is called a flag of weak maximal contact for $I_X$ at $x$ if it
maximizes the resolution invariant lexicographically among all choices of 
admissible flags at $x$.
\end{df}

This definition obviously behaves well under passage to a coefficient ideal
w.r.t. $H_1$ by omitting the first entry $H_1$ from ${\mathcal H}$ to obtain
the new flag ${\mathcal H}_{H_1}$. This is again a flag of maximal contact,
since conditions (a)-(c) and maximality follow trivially from the respective 
conditions on ${\mathcal H}$. Hence considering a flag of weak maximal contact 
instead of a hypersurface of weak maximal contact does not change any of the 
key properties, but allows more flexibility for dealing with lower level
kangaroos. 

\section{Two different kinds of double kangaroos}

It is a well known fact that the situation in positive characteristic 
can only differ from the one in characteristic zero in rather special
situations. 
Hauser studied such phenomena in great detail in \cite{Ha} 
by considering precisely the two levels involved in a kangaroo point.
For surfaces, he and Wagner extended these considerations to a 
general treatment of the purely inseparable case in \cite{HW}.
The situation in higher dimension differs from this easiest case
in the sense that there might be more than just two levels at which 
the ridge is not generated in degree 1 at some time during the process of
blowing ups. The following two examples illustrate three different 
roles of the different levels of the flag of weak maximal contact in such a 
setting.

\begin{df}
Let ${\mathcal H}$ be a flag of weak maximal contact for an ideal 
$I_X \subset W$ at the point $x$ which we assume for simplicity to
be the origin of our coordinate chart. We denote the $i$-th coefficient 
ideal, which arises when descending to $H_i$, by 
$J_i \subset {\mathcal O}_{H_i}$. If the ideal generated by the lowest
order generators of $J_{i-1}$ is not a principal ideal,
$H_i$ is called 
\begin{itemize}
\item {\bf neutral}, if the degree $1$ part of the generator of the principal 
      ideal $I_{H_i} \subset {\mathcal O}_{H_{i-1},0}$ is in the 
      ${\mathbb C}$-span of the degree $1$ elements of the ridge/n-ridge
      of $J_{i-1}$.
\item {\bf active}, if it is the $H_i$ of lowest index $i$ which is not
      neutral.
\item {\bf dormant}, if it is neither active nor neutral.
\end{itemize}
If, on the other hand, the ideal generated by the lowest order generators
of $J_{i-1}$ is principal, it is of the form $g^{\frac{b!}{b-k}}$ for some 
$k<b$ and we change the notions of neutral, active and dormant by replacing
the ridge/n-ridge of $J_{i-1}$ by the one of $\langle g \rangle$. 
\end{df}

\begin{rem}
\begin{enumerate}
\item According to Hauser's description of the process leading to kangaroo 
      points, at least one active $H_i$ and one dormant $H_j$ are necessary 
      to produce a kangaroo phenomenon.
\item If the ideal generated by the lowest order generators of $J_{i-1}$ is 
      not principal, there is at least one ideal among the
      contributing $I_k$, of which the ideal generated by its lowest order
      generators is itself not principal, e.g. generated by $f_1$ and $f_2$.
      Hence taking the $\frac{b!}{b-k}$-th power of of this $I_k$ upon 
      forming the coefficient ideal, we obtain all mixed products of the 
      form $f_1^af_2^b$, $a+b=\frac{b!}{b-k}$. This implies that higher 
      degree generators of the n-ridge can only occur if they would also
      occur for $\langle f_1, f_2 \rangle$.\\
      If on the other hand, the ideal generated by the lowest order generators
      of $J_{i-1}$ is principal, the generator is of the form 
      $g^{\frac{b!}{b-k}}$
      for some $k$ and hence masks the true situation of the (n-)ridge of $g$.
      This is the reason for the special treatment of this case in the
      above definition.
\end{enumerate}
\end{rem}
 
Both of the following examples were constructed in a straight forward way, 
combining two occurrences of kangaroos at two different levels.  
Similar examples can be constructed in any positive characteristic 
and for any ambient dimension exceeding 4. However, these examples involve 
several blow-ups between the first and the second occurrence, basically 
making a fresh start after the first. Here no effort is made to reduce 
this number of blow-ups, since the context of this article is the study 
of the roles of the hypersurfaces of weak maximal contact. \\

To keep these rather lengthy examples more readable, we only state the
blow-ups, the weak transform at each step and the flag of weak maximal
contact, whenever the latter changes, but omit all data
which is related to coefficient ideals, since these can easily be 
computed for these examples.

\begin{example} \label{1312}
In this example, a hypersurface in ${\mathbb A}_K^5$, $char K=3$,
we shall see 2 occurrences of kangaroo points on two different levels of 
coefficient ideals. For both occurrences, the active hypersurface of weak 
maximal contact is the first one in the flag. Note that the two blowing ups 
with chart $E=V(y)$ after the first kangaroo are only used for setting up 
the degrees for the following kangaroo.\footnote{Whenever we write 'h.o.t.'
we want to indicate that there are further terms of higher degree, which are
irrelevant for the further considerations. In this case only the first 
non-relevant term is stated, even if this does not happen to be the term 
originating from the previous first non-relevant term}\\

\begin{itemize}
\item before 1st blowing up\\
      $I=\langle w^3+y^6z^3v^2+x^9y^8+x^{18}y^2+x^{18}v^2 \rangle$\\
      Flag: \\
      $V(w)$ active, $V(w,v)$ neutral, $V(w,v,z)$ dormant, 
            $V(w,v,z,y)$ neutral
\item after blowing up at the origin, chart $E_{new}=V(x)$\\
      $I=\langle w^3+x^8(y^6z^3v^2+x^6y^8+x^9y^2+x^9v^2) \rangle$
\item after blowing up at the origin, chart $E_{new}=V(y)$\\
      $I=\langle w^3+x^8y^{16}(z^3v^2+x^6(x^3+y^3)+x^9v^2) \rangle$\\
      Flag: $V(w)$ active, $V(w,v)$ neutral, $V(w,v,z)$ dormant, $V(w,v,z,y)$ dormant
\item after blowing up at the origin, chart $E_{new}=V(x)$\\
      $I=\langle w^3+x^{26}y^{16}(z^3v^2+x^4+x^4y^3+x^6v^2) \rangle$\\
      coordinate change: $y_{new}=y_{old}+1$, $w_{new}=w_{old}+x^{10}y$\\
      $I=\langle w^3+x^{26}((y-1)^{16}(z^3v^2+x^6v^2)-x^4y^4 + h.o.t.)
         \rangle $\\
      Flag in new coordinates: \\
      $V(w)$ active, $V(w,v)$ neutral, $V(w,v,z)$ dormant, $V(w,v,z,y)$ neutral\\
      Kangaroo at 3rd coefficient ideal
\item after blowing up at the origin, chart $E_{new}=V(x)$\\
      $I=\langle w^3+x^{28}((xy-1)^{16}(z^3v^2+x^3v^2)-x^3y^4 +h.o.t.)
         \rangle$
\item after blowing up at the origin, chart $E_{new}=V(v)$\\
      $I=\langle w^3+x^{28}v^{30}(z^3+x^3-x^4yv^2 +h.o.t.)\rangle$
\item after blowing up at the origin, chart $E_{new}=V(y)$\\
      $I=\langle w^3+x^{28}y^{58}v^{30}(x^3+z^3-x^4y^4v^2+h.o.t.)\rangle$
\item after blowing up at the origin, chart $E_{new}=V(y)$\\
      $I=\langle w^3+x^{28}y^{116}v^{30}(x^3+z^3-x^4y^7v^2+h.o.t.)\rangle$
\item after blowing up at the origin, chart $E_{new}=V(z)$\\
      $I=\langle w^3+x^{28}y^{116}z^{174}v^{30}(1+x^3-xy^7z^{10}v^2+h.o.t.)
         \rangle$\\
      coord. change: $x_{new}=x_{old}+1$, $y_{new}=y_{old}+1$,
                         $w_{new}=w_{old}+z^{58}v^{10}x^3$\\
      $I=\langle w^3+z^{174}v^{30}(x^4+x^3y+h.o.t) \rangle$\\
      Flag: $V(w)$, $V(w,x)$, $\dots$\\
      Kangaroo at 1st coefficient ideal
\end{itemize}
\end{example}

\begin{example} \label{2312}
In this example, again in the same affine space as before, we shall see 2 
occurrences of kangaroo points on two different levels of coefficient ideals. 
For the first occurrence, a dormant hypersurface of weak maximal contact 
acts as the active one, for the second it is the top-level active 
hypersurface of weak maximal contact. 
This example again basically consists of two regular kangaroo 
phenomena in a row, occurring on two different levels, but in a different
flavor than example \ref{2312}.\\

\begin{itemize}
\item before first blowing up\\
      $I=\langle w^3+xy^9z^9v+x^7y^{20}v+x^{34}y^2v+x^{46}v \rangle$\\
      Flag:\\
      $V(w)$ active, $V(w,v)$ neutral, $V(w,v,z)$ dormant, $V(w,v,y,z)$ neutral
\item after blowing up at the origin, chart $E_{new}=V(x)$\\
      $I=\langle w^3+x^{17}(y^9z^9v+x^8y^{20}v+x^{17}y^2v+x^{27}v) \rangle$
\item after blowing up at the origin, chart $E_{new}=V(y)$\\
      $I=\langle w^3+x^{17}y^{33}(z^9v+x^8y^{10}v+x^{17}yv+x^{27}y^9v) 
         \rangle$\\
      Flag: $V(w)$ active, $V(w,v)$ neutral, $V(w,v,z)$ dormant, $V(w,v,y,z)$ dormant
\item after blowing up at the origin, chart $E_{new}=V(x)$\\
      $I=\langle w^3+x^{57}y^{33}(z^9v+x^9y^{10}v+x^9yv+x^{27}y^9v) 
         \rangle$\\
      coord. change: $y_{new}=y_{old}+1$, $z_{new}=z_{old}+xy$\\
      $I=\langle w^3+x^{57}(y-1)^{33}(z^9v+x^9y^{10}v-x^{27}v+x^{27}y^9v)
         \rangle$\\
      Flag in new coordinates: \\
      $V(w)$ active, $V(w,v)$ neutral, $V(w,v,z)$ dormant, 
      $V(w,v,y,z)$ neutral\\
      Kangaroo at 3rd coefficient ideal
\item after blowing up at the origin, chart $E_{new}=V(x)$\\
      $I=\langle w^3+x^{64}(xy-1)^{33}(z^9v+x^{10}y^{10}v-x^{18}v+x^{27}y^9v)
         \rangle$
\item after blowing up at the origin, chart $E_{new}=V(x)$\\
      $I=\langle w^3+x^{71}(x^2y-1)^{33}(z^9v+x^{11}y^{10}v-x^9v+x^{27}y^9v)
         \rangle$
\item after blowing up at the origin, chart $E_{new}=V(v)$\\
      $I=\langle w^3+x^{71}(x^2yv^3-1)^{33}v^{78}(z^9-x^9 +h.o.t) \rangle$\\
      Flag: $V(w)$ active, $V(w,z)$ dormant, $V(w,x,z)$ dormant, $\dots$
\item after blowing up at the origin, chart $E_{new}=V(y)$\\
      $I=\langle w^3+x^{71}y^{155}v^{78}(x^2y^6v^3-1)^{33}(z^9-x^9+h.o.t)
         \rangle$
\item after blowing up at the origin, chart $E_{new}=V(y)$\\
      $I=\langle w^3+x^{71}y^{310}v^{78}(x^2y^6v^3-1)^{33}(z^9-x^9+h.o.t)
         \rangle$
\item after blowing up at the origin, chart $E_{new}=V(z)$\\
      $I=\langle w^3+x^{71}y^{310}v^{78}z^{465} (x^2y^6v^3z^{16}-1)^{33}
                 (1-x^9+h.o.t.) \rangle$\\
      coord. change: $x_{new}=x_{old}-1$, $y_{new}=y_{old}+1$, 
                     $w_{new}=w_{old}-v^{26}z^{155}x^3$ \\
      Kangaroo at 1st coefficient ideal
\end{itemize}

\end{example}

In both examples the relevant order of the first respectively second
coefficient ideal dropped significantly after the first kangaroo 
phenomenon, but before the occurrence of the kangaroo on this level. 
The examples have been constructed to illustrate roles of hypersurfaces 
of maximal contact in multiple kangaroos and not to specifically illustrate 
the increase in order. Nevertheless the observed behaviour raises several 
questions, which seem to be natural starting points for further experiments
in the search for new meaningful examples:\\
\begin{itemize}
\item Is it possible to find an occurrence of two kangaroo phenomena
      whose 'distance' is less than 3 blow ups?
\item Is it possible to find an occurrence of two kangaroo phenomena
      for which the drop of order between the first and the second kangaroo 
      does not outweigh the increase of order?
\item If one of the previous question has an affirmative answer, what is
      the smallest dimension in which this occurs?
\end{itemize}

\end{document}